# The Story of Hagge and Speckman

**Christopher J Bradley**

**Abstract:** This concerns theory first developed by Hagge and Speckman in the Edwardian era. Speckman investigated triangles that were simultaneously in perspective and indirectly similar. On the other hand Hagge studied circles that pass through the orthocentre of a given triangle. Superficially these subjects look unrelated, but this is not the case, as we describe.

## 1. Introduction

In 1905 Speckman [1] wrote a paper on indirectly similar triangles in perspective and in1907 Hagge [2] wrote another paper on the construction of circles that always pass through the orthocentre of a given triangle, which have become known as Hagge circles. It is evident that Hagge was unaware of Speckman's work, and indeed it appears that the contents of Speckman's paper, possibly because it is written in Dutch, seem to be little known. Hagge's work, on the other hand has been developed, notably by Peiser [3], and more lately by Bradley and Smith [4], [5]. Both the 1905 and 1907 papers are ground breaking papers in classical geometry, and the work of Speckman, in particular, deserves to be better known. The interesting thing is that the two pieces of work are related in a manner that at first sight is not at all obvious. However, if Hagge had read Speckman's paper it is almost certain that he would have understood its relevance and we would not have had to wait so long to see the developments of Hagge's work that are a consequence. For reasons that will become clear it is best to describe Hagge's work first and then to summarise Speckman's work, showing how it may be applied, step by step, to Hagge's configuration.

## 2. Hagge circles

For the construction of a Hagge circle we refer to Fig.1.

Given a triangle *ABC* and a point *P* not on a side or the extension of a side and not on the circumcircle$\Gamma$, the Hagge circle of *P* with respect to *ABC* is defined as follows: draw *AP*, *BP*, *CP* to meet $\Gamma$ at *D*, *E*, *F* respectively; now reflect *D*, *E*, *F* in the sides *BC*, *CA*, *AB* respectively to obtain the points *U*, *V*, *W*. The Hagge circle $\Sigma(P)$ is now defined to be the circle *UVW*. In his paper [2] Hagge proved

(i) that $\Sigma(P)$ always passes through *H*, the orthocentre of triangle *ABC*;
(ii) that, if $\Sigma(P)$ intersects the altitudes *AH*, *BH*, *CH* at *X*, *Y*, *Z* respectively, then *UPX*, *VPY*, *WPZ* are straight lines;
(iii) that the Hagge circle of *G*, the centroid of *ABC*, with respect to the medial triangle (the triangle joining the midpoints of the sides) is the Brocard circle on *OK* as diameter, where $K$ is the symmedian point of $ABC$, $O$ being the circumcentre of $ABC$ and thus the orthocentre of the medial triangle. At the time, Brocard geometry was a subject of particular interest, which may account for why roughly half of Hagge's paper is devoted to this application.

Peiser's [3] main contribution was concerned with the centre $Q$ of the circle $\Sigma(P)$. Using the algebra of complex numbers he was able to show that the image of $Q$ under a half-turn about the nine-point centre $T$ is the same point as the isogonal conjugate $Pg$ of $P$ with respect to triangle $ABC$. Since $T$ is the midpoint of $OH$ it follows that the figure $PgHQO$ is a parallelogram. Another property not mentioned by earlier geometers is that the midpoints of $AU$, $BV$, $CW$ all lie on the nine-point circle.

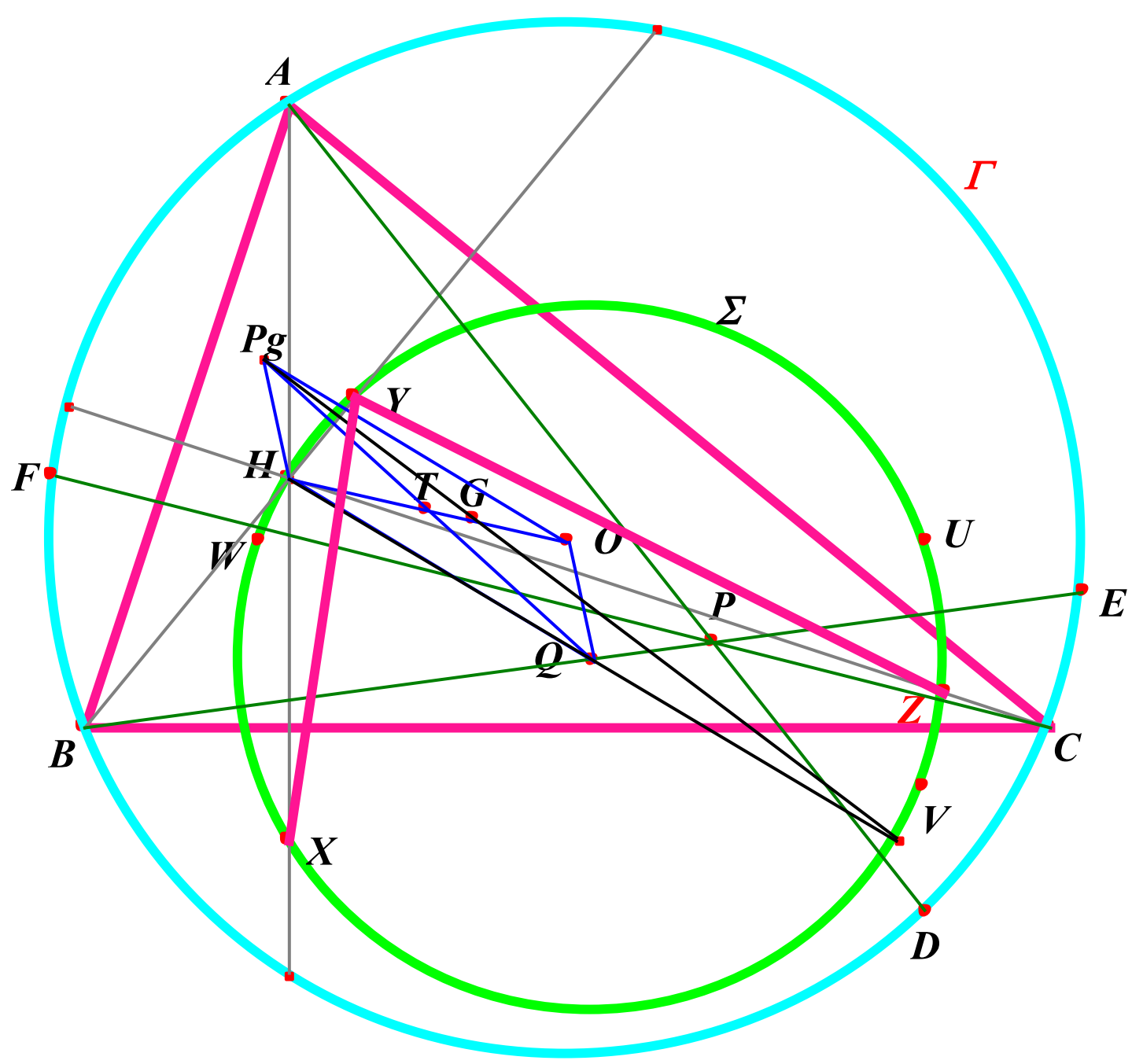

**Fig. 1**
**The construction of a Hagge circle**

If $P$ lies on the circumcircle the lines $APg$, $BPg$, $CPg$ are parallel, so $Pg$ lies on the line at infinity. The consequence is that the Hagge circle degenerates into a line through $H$, which being the line through the reflections of $P$ in the sides $BC$, $CA$, $AB$ is parallel to the Simson line of $P$ and is sometimes called the double Simson line of $P$.

When $P$ is the incentre $I$ the Hagge circle is the *Fuhrmann* circle and this circle was known about many years before Hagge's publication. The Fuhrmann circle not only passes through $H$ but also passes through Nagel's point $Na$. An account of the Fuhrmann circle is given by Honsberger [5], who shows that $HNa$ is one of its diameters. When $P$ is the symmedian point $K$ the Hagge circle is the orthocentroidal circle on $GH$ as diameter. In general there is a straightforward result that *PgG* meets the Hagge circle at the point at the other end of the diameter through *H*. See Fig. 1.

In their first paper Bradley and Smith [4] give a synthetic proof of Peiser's result and also show how a consistent definition can be given for a Hagge circle if *P* lies on a side of the triangle or at a vertex. In their second paper Bradley and Smith [5] provide another proof of Hagge's result (ii).

Now if you consider the quadrilateral *HXYZ*, since *CZH* and *BHY* are perpendicular to the sides, we have $\angle ZXY = \angle ZHY = \angle BAC$. Similar angle relations show that triangles *ABC* and *XYZ* are similar, and since the ordering of the letters in the two triangles is opposite, these triangles are indirectly similar. The same applies to triangles *UVW* and *DEF*. Bradley and Smith [4] also give a synthetic proof of the important result that there is an indirect similarity carrying *ABCDEFP* to *XYZUVWP* in which *P* is the unique fixed point. In Fig. 2 we provide an illustration of the dilative reflection through *P* and its axis that provides this indirect similarity.

On reflection in the axis through *P* the figure *ABCDEF* is mapped to $X'Y'Z'U'V'W'$ and then by dilation through the centre *P* by a factor *XY*/*AB* the figure $X'Y'Z'U'V'W'$ is mapped to *XYZUVW*. Note that if the perpendicular axis is used instead, then a rotation of $180^{o}$ is required in addition to the dilation and reflection.

It is also significant, as we shall see in due course, that triangles *ABC* and *XYZ* are in perspective with vertex of perspective *H*, and that triangles *XYZ* and *UVW* are in perspective with vertex of perspective *P*.

Further results proved by Bradley and Smith [5] are

(i) that, if *VW* meets *AH* at $U'$, *WU* meets *BH* at $V'$ and *UV* meets *CH* at $W'$, then $U'$, $V'$, $W'$, *P* are collinear, and

(ii) that, if $U_1$, $V_1$, $W_1$, $X_1$, $Y_1$, $Z_1$ are the midpoints of *DU*, *EV*, *FW*, *AX*, *BY*, *CZ*, then $U_1$, $V_1$, $W_1$, $X_1$, $Y_1$, $Z_1$ lie on a conic, which they call the *midpoint conic*.

Also the midpoints of *AX*, *BY*, *CZ*, *DU*, *EV*, *FW* lie on a conic. In fact one can go further and say that the points that divide these six segments in any fixed ratio lie on a conic.

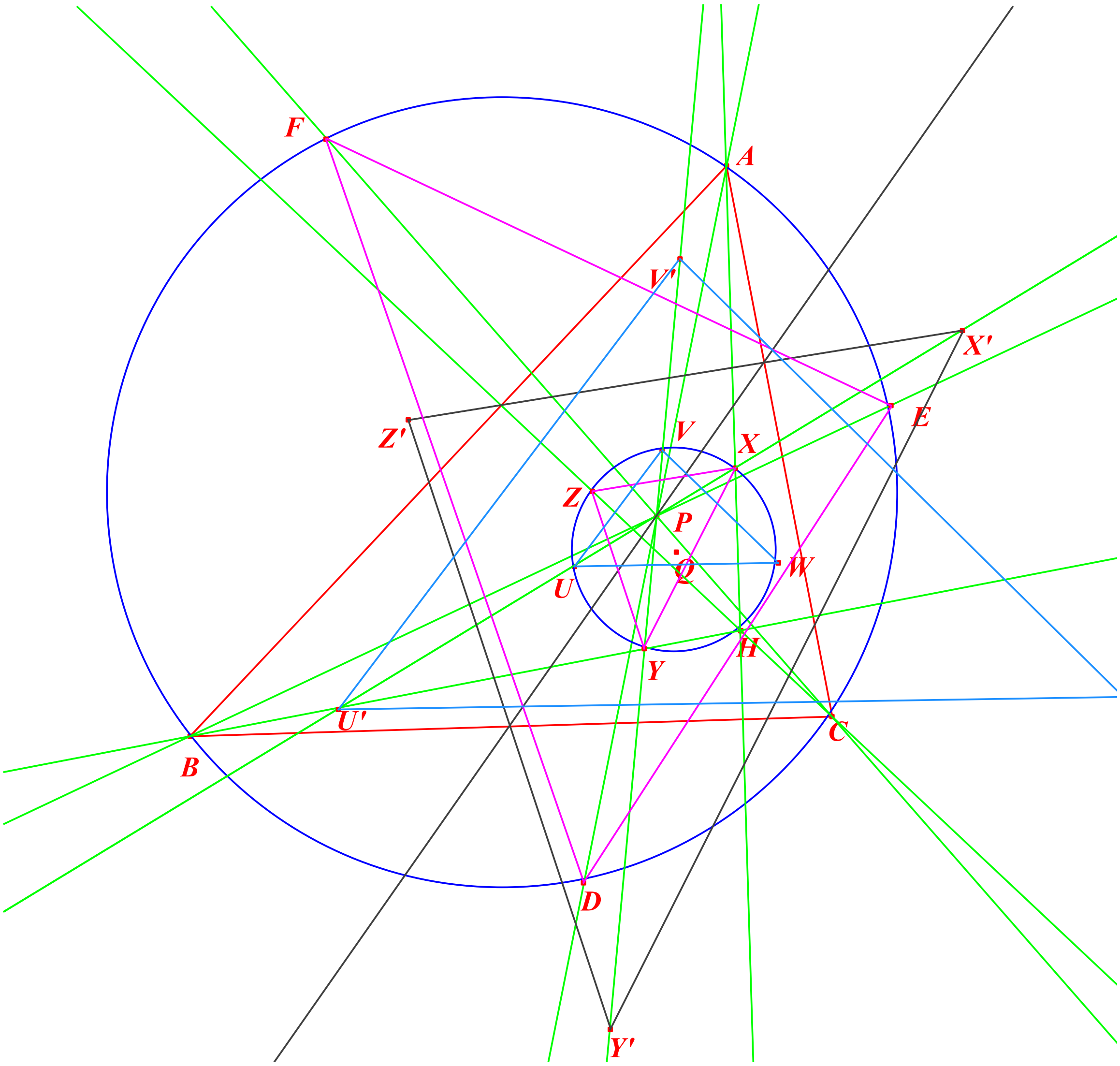

**Fig. 2**
**Some properties of a Hagge circle**

## 3. Speckman's theorems on indirectly similar triangles in perspective

Speckman announced proofs of his theorems in April 1903 at a Congress in ś-Gravenhage and they were published [1] in 1905. Their relevance to the work of Hagge lies in the fact that in the Hagge configuration, described in Section 2, triangles *ABC* and *XYZ* are indirectly similar and are in perspective with vertex of perspective *H*.

If Hagge had read Speckman's work, then it seems highly likely that he would have applied the theorems to obtain a sequence of additional results. And if Speckman had known about Hagge's work it seems almost certain he would have referred to it subsequently. As far as we are aware the application of Speckman's results to the Hagge configuration has never been pointed out and is a story waiting to be told. We tell this story by relating the more important of Speckman's results paragraph by paragraph, pointing out in each case its significance in relation to Hagge's configuration. Whilst it is true in an obvious sense that Speckman's results are more general than Hagge's (he deals with the case when the vertex of perspective is a general point and not the orthocentre and consequently derives more general results), it is possible to hold an alternative point of view that Hagge's configuration is illustrative of all the properties one might need to know about pairs of indirectly similar triangles in perspective. For if you have an arbitrary pair of indirectly similar triangles in perspective, one of them is always related to the other by being directly similar to one of its partner's *Hagge triangles* (defined to be the triangle *XYZ* inscribed in a Hagge circle).

This is illustrated in Fig. 3. It shows two arbitrary triangles *ABC* and *A′B′C′* that are indirectly similar by means of a reflection in the line *l* through *P* followed by dilation through *P*. They are chosen to be in perspective with vertex a point *Q*. Triangle *XYZ* is the Hagge triangle of *P* and triangle *X′Y′Z′* shows a dilation of *XYZ* through *P* with an appropriate scale factor. It can now be observed that triangles *A′B′C′* and *X′Y′Z′* are images of one another under a rotation about *P*. If the point *P* is not known it may be obtained as the other point of intersection besides *Q* of the two rectangular hyperbolae *ABCHQ* and *A′B′C′H′Q*. See the first paragraph of Section 4.

## 4. Paragraphs 1 to 9 of Speckman's paper

We refer to Fig. 3 in which *H* and *H′* are the orthocentres of triangles *ABC* and *A′B′C′*, *Q* is the perspector and *P* is what Speckman calls the double-point of inverse similarity. He proves in the first paragraph that rectangular hyperbolas may be drawn passing through *ABCHPQ* and *A′B′C′H′PQ* and possessing asymptotes that are parallel.

From now on we describe Speckman's results primarily as they apply to the Hagge configuration. The more general results, if not stated, should be clear enough.

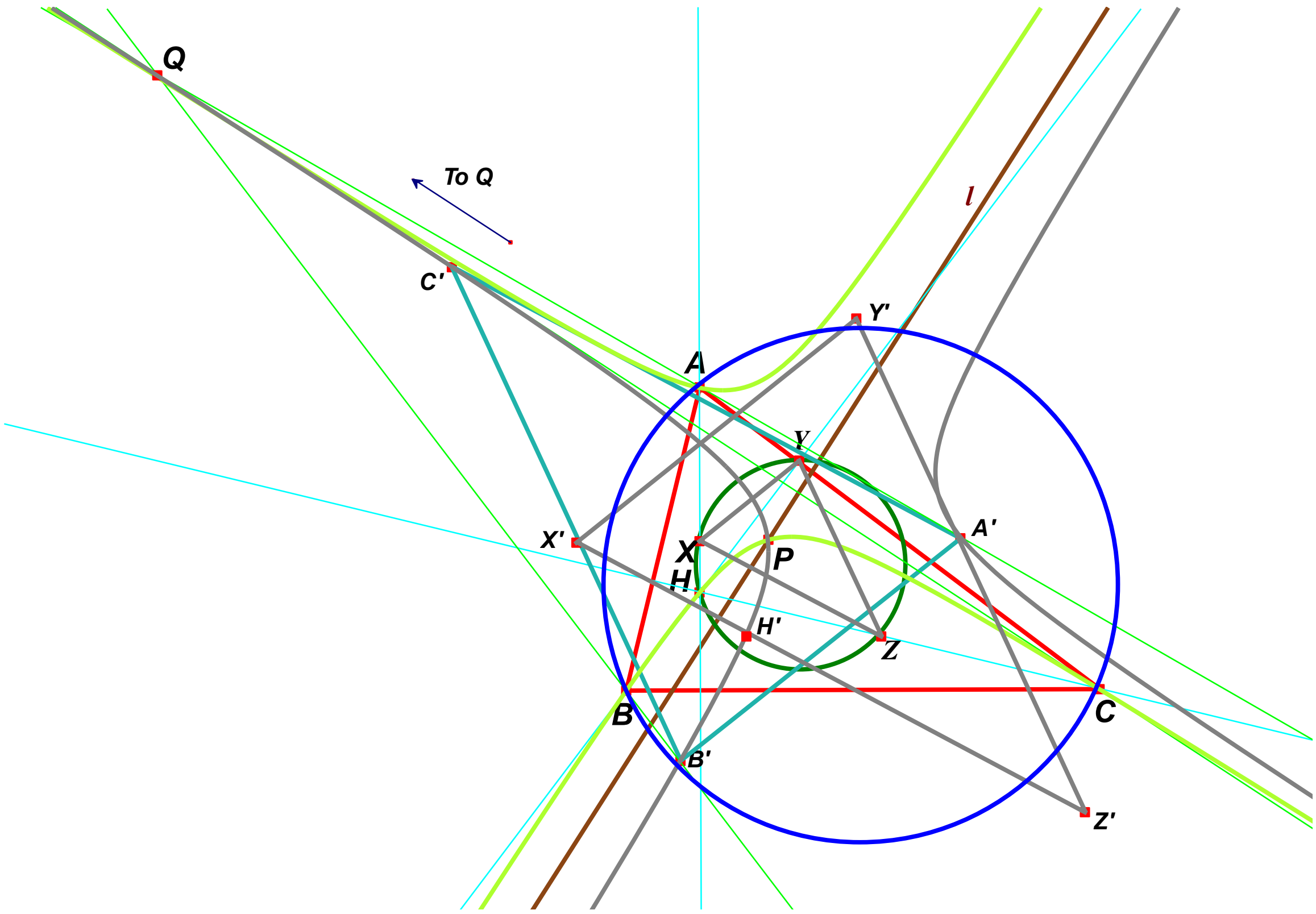

**Fig. 3**
**Indirectly similar in perspective triangles are related to Hagge triangles**

In Fig. 4 the result of the first paragraph is illustrated when the second triangle is the Hagge triangle *XYZ*. Note that *Q* now coincides with *H* and the asymptotes are parallel to the axes through *P*, which Speckman calls the *double lines of inverse similarity*. Also he shows in the second paragraph that these lines are parallel to the angle bisectors of corresponding lines of triangles *ABC* and *XYZ*. The reflection of *XYZ* in one of these axes is the triangle *X′Y′Z′*, which is a dilation of triangle *ABC* through *P*. Triangles *ABC* and *XYZ* are orthologic, as indeed are all pairs of indirectly similar triangles (a fact we prove later in this article) and the orthologic centres of the triangles are *D* and *H*, where the point *D* lies on both the circumcircle of *ABC* and the rectangular hyperbola *ABCHP*. Also, as shown in Bradley and Smith [4], the orthocentre *h* of triangle *XYZ*, *P* and *D* are collinear. The centres of the two rectangular hyperbolas are denoted by *M* and *m*. Also illustrated in Fig. 4 is the fact, proved in Speckman's third paragraph, that *J*, *j′* and *P* are collinear where *j′* is the $180^{\circ}$ rotation of *j* about *m*, and *J*, *j* are any pair of corresponding points on the two hyperbolae under the similarity.

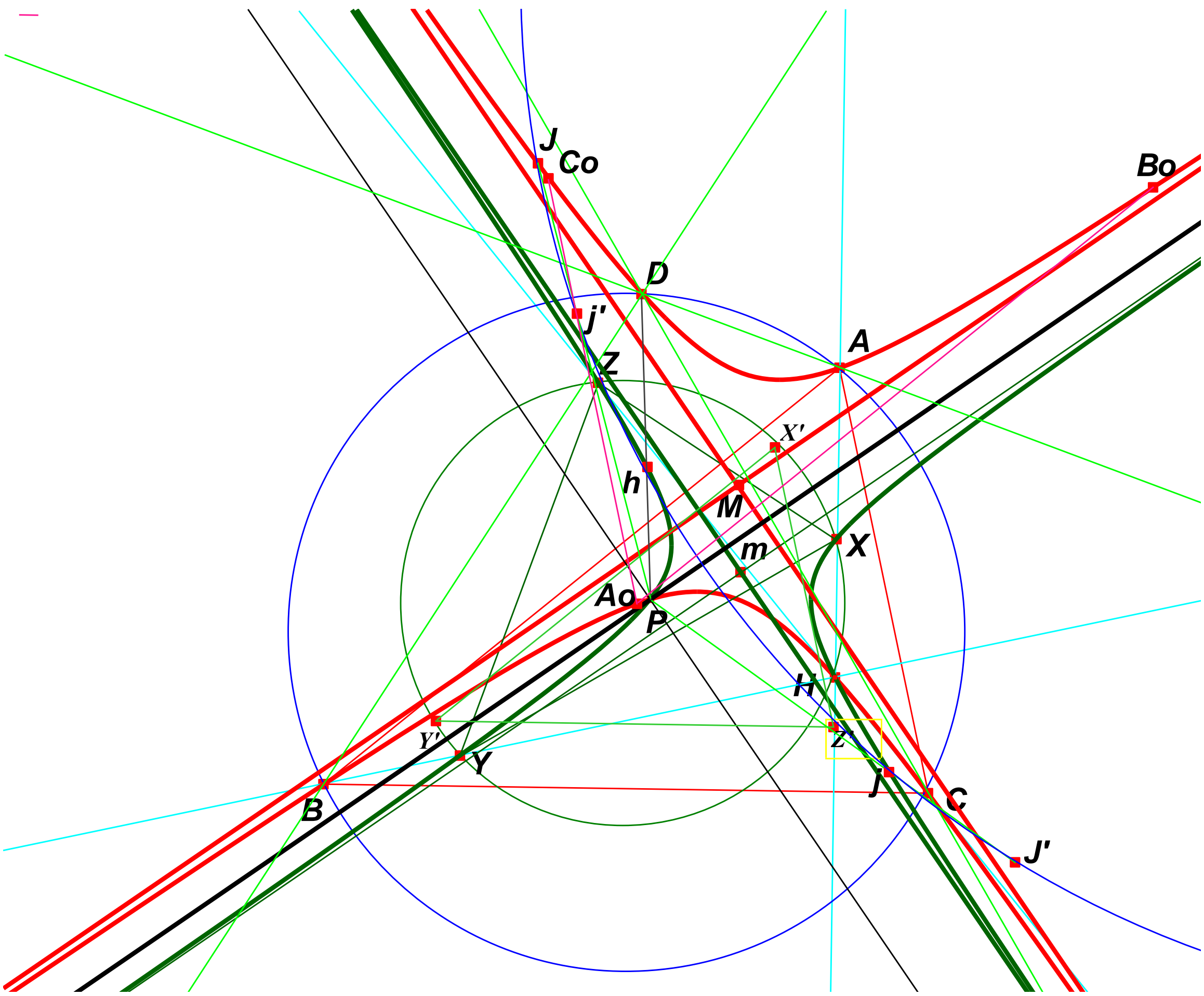

**Fig. 4**
**Illustrating the first three paragraphs of Speckman's paper**

In paragraph four there is a construction, whereby triangle *ABC* is rotated by 180° about the centre *M* of the hyperbola *ABCHP* to produce a triangle $A_0B_0C_0$ and it is proved that triangles *XYZ* and $A_0B_0C_0$ are also inversely similar and in perspective, but now the vertex of perspective is *P*, the double point of inverse similarity is *H* and the orthocentre of triangle $A_0B_0C_0$ is *D*. In the more general case the interchange of the roles of the double point of inverse similarity *P* and the vertex of perspective *Q* arising from this construction is an intriguing result.

In paragraph five there is proof of the result that corresponding points on the two hyperbolae together with the points obtained by half-turns about their respective centres *M* and *m* are concyclic. An arc of the circle *Jj′jJ′* can be seen in Fig. 4.

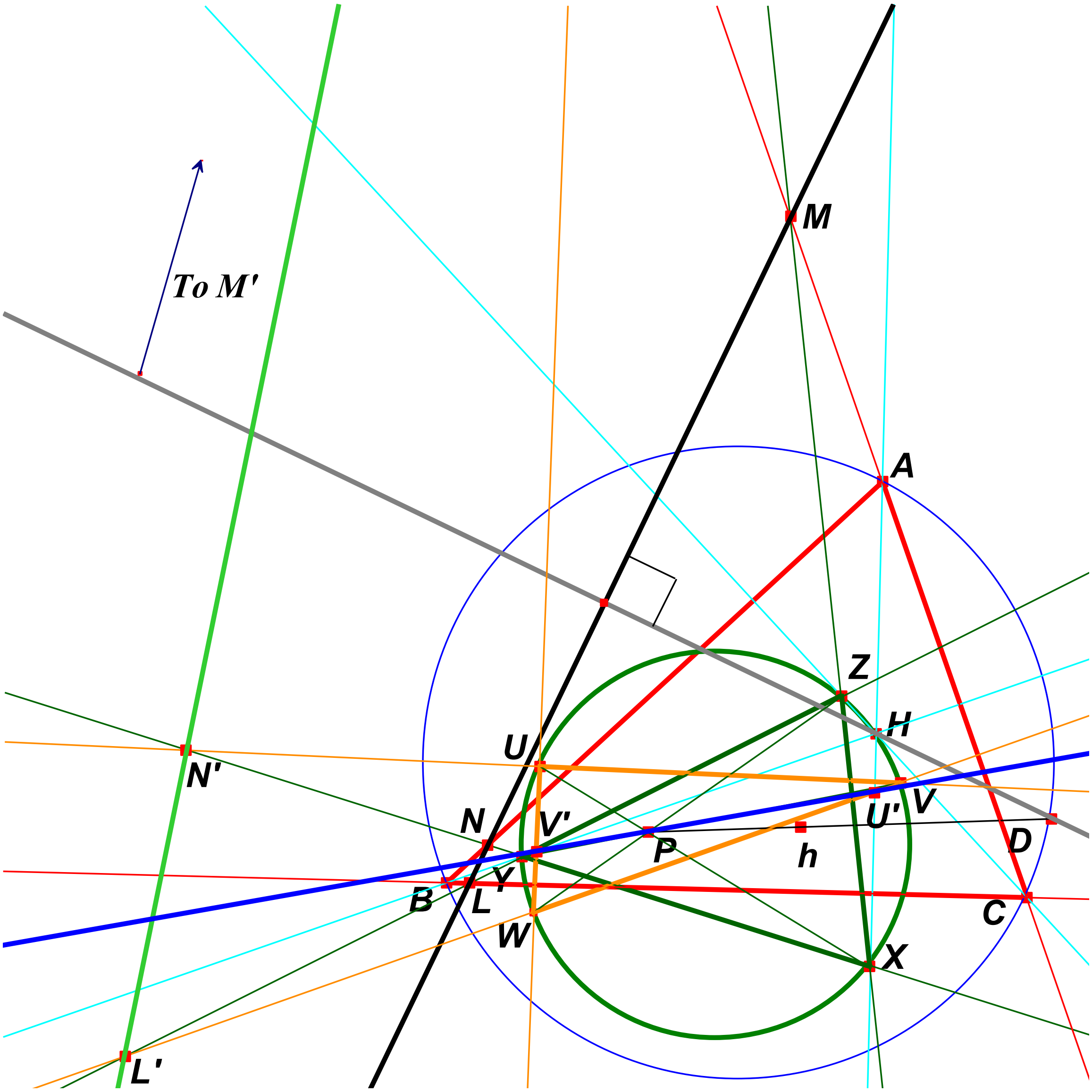

**Fig. 5**
**Desargues' axes of perspective and the orthological centres**

In paragraph seven there is a proof that in the general case (when the perspector is not at one of the orthocentres) the two orthocentres are collinear with the perspector. This result is trivial in the Hagge configuration, since $Q$ then coincides with $H$. In paragraph nine it is proved that the lines connecting the orthological centre of one triangle to the orthocentre of the other passes through the double point of inverse similarity. When applied to the Hagge configuration this implies that the line $hD$ passes through $P$, a result which is also proved in Bradley and Smith [5].

## 5. Axes of perspective in the Hagge configuration

In Fig. 5 we show the axis *U′V′W′* defined at the end of Section 2 and also the Desargues axes *LMN*, *L′M′N′* of perspective arising respectively from the perspectives of triangles *ABC*, *XYZ* and of triangles *XYZ*, *UVW*. Paragraph ten of contains the proof that the axis *LMN* is perpendicular to the line *HD* joining the orthologic centres.

## 6. Paragraphs 11 to 18 of Speckman's paper

Fig. 6 illustrates the results of paragraphs 11-16. Most of the notation has previously been defined. *Pg* is the isogonal conjugate of *P* with respect to triangle *ABC* and *pg* is the isogonal conjugate of *P* with respect to triangle *XYZ*.

In paragraph eleven it is noted that the line through *A* parallel to *YZ*, the line through *B* parallel to *ZX* and the line through *C* parallel to *XY* are known to meet and it is proved they do so at a point *S* on the circumcircle of *ABC*. Speckman calls this point the *metapole* of triangle *ABC* with respect to triangle *XYZ*. Such a point is now called a *paralogic centre*, so it turns out that we are dealing with a case in which the triangles are not only orthologic but also paralogic. (It is not always the case that orthologic triangles are also paralogic.)The point *s* is defined similarly. In the next two paragraphs there are proofs about the properties of the lines *DS* and *Hs*. They are diameters, respectively, of circles *ABC* and *XYZ*. It is proved that the line *Hs* is parallel to the isogonal conjugate of the rectangular hyperbola *ABCHP*, which is a line through *Pg*, and that *DS* is parallel to the isogonal conjugate of the rectangular hyperbola *XYZHPh*, which is a line through *pg*. Speckman calls the lines *DS* and *Hs* orthological diameters of the circles. A curious fact is that *in the Hagge configuration* the intersection *R* of the two orthological diameters lies on the rectangular hyperbola *ABCHP*.

In paragraph fourteen the midpoints of the sides of *ABC* are denoted by $M_a$, $M_b$, $M_c$ so that triangle $M_aM_bM_c$ is indirectly similar to triangle *XYZ*. It is shown that the perpendiculars from $M_a$ on to *YZ*, $M_b$ on to *ZX* and from $M_c$ on to *XY* meet at a point *W*, which is the midpoint of *HS*. The point *w* is similarly defined, and is the midpoint of *hs*. Also $HS/hs = AB/XY$.

In paragraph fifteen it is proved that the paralogic centre of triangle $M_aM_bM_c$ with respect to triangle *XYZ* is the midpoint of *DH*. This is the centre of the hyperbola *ABCHPD*. In paragraph seventeen it is shown that the Desargues axis of perspective *LMN*, see Fig. 5, bisects the line segment *Ss* joining the paralogic centres.

In paragraph sixteen it is shown that if you reflect triangle *XYZ* in the Desargues axis *LMN* to get triangle *X′Y′Z′*, then triangle *X′Y′Z′* is in perspective with triangle *ABC* with vertex of perspective the point where *DH* meets the circumcircle *ABC*. In paragraph eighteen it is then noted that if the axis *LMN* meets the common chord of circles *ABC* and *X′Y′Z′* at *E*, then *E* is the radical centre of circles *ABC*, *XYZ* and *X′Y′Z′*. Hence the line through *E* perpendicular to *OQ* is the radical axis of the circumcircle and the Hagge circle. The construction of the radical axis of the circumcircle and the Hagge circle is a splendid flourish with which to end the paper.

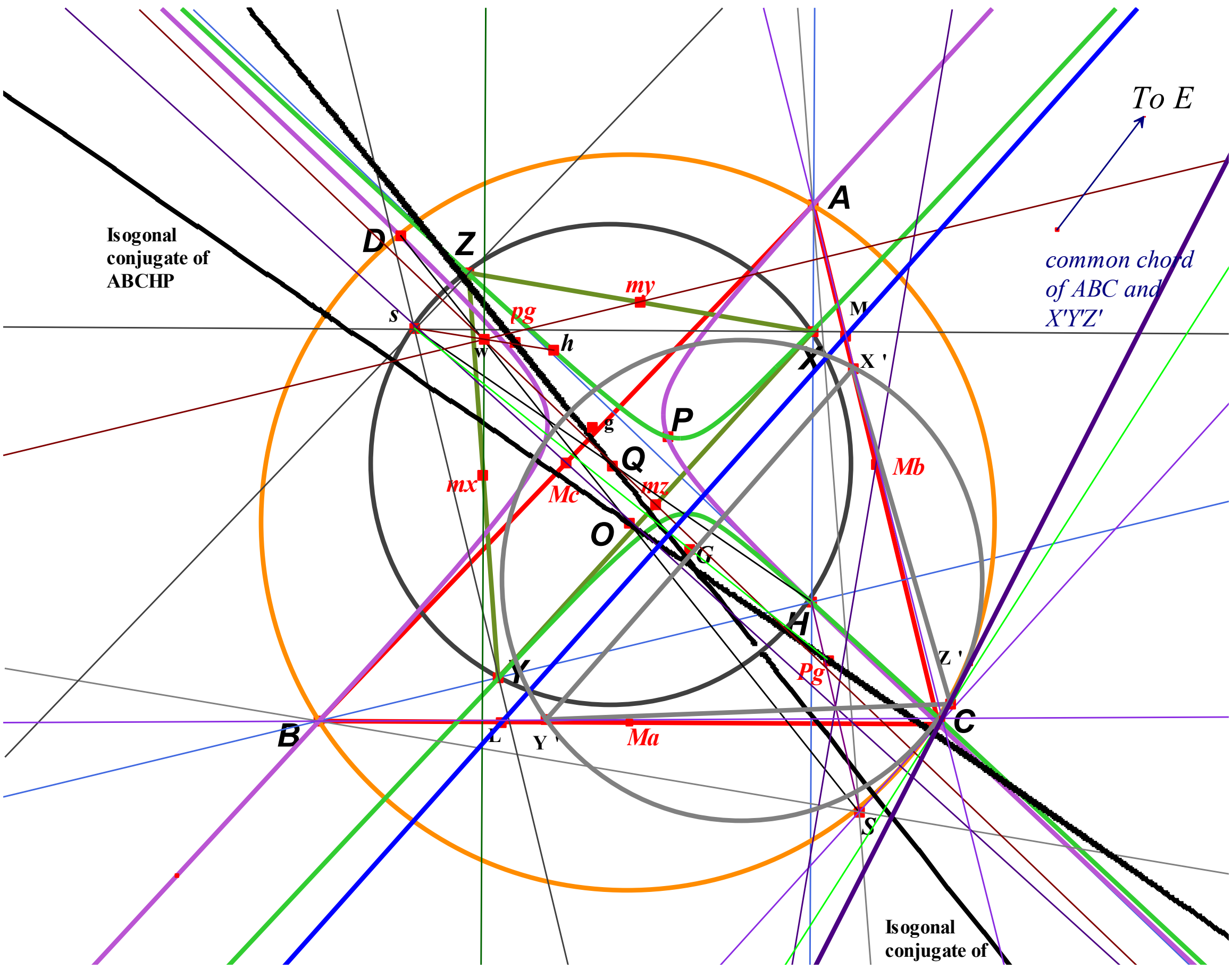

**Fig.6**
**Illustrating paragraphs 11-16 of Speckman's paper**

7. **Which of all these properties hold if there are two indirectly similar triangles that are not necessarily in perspective?**

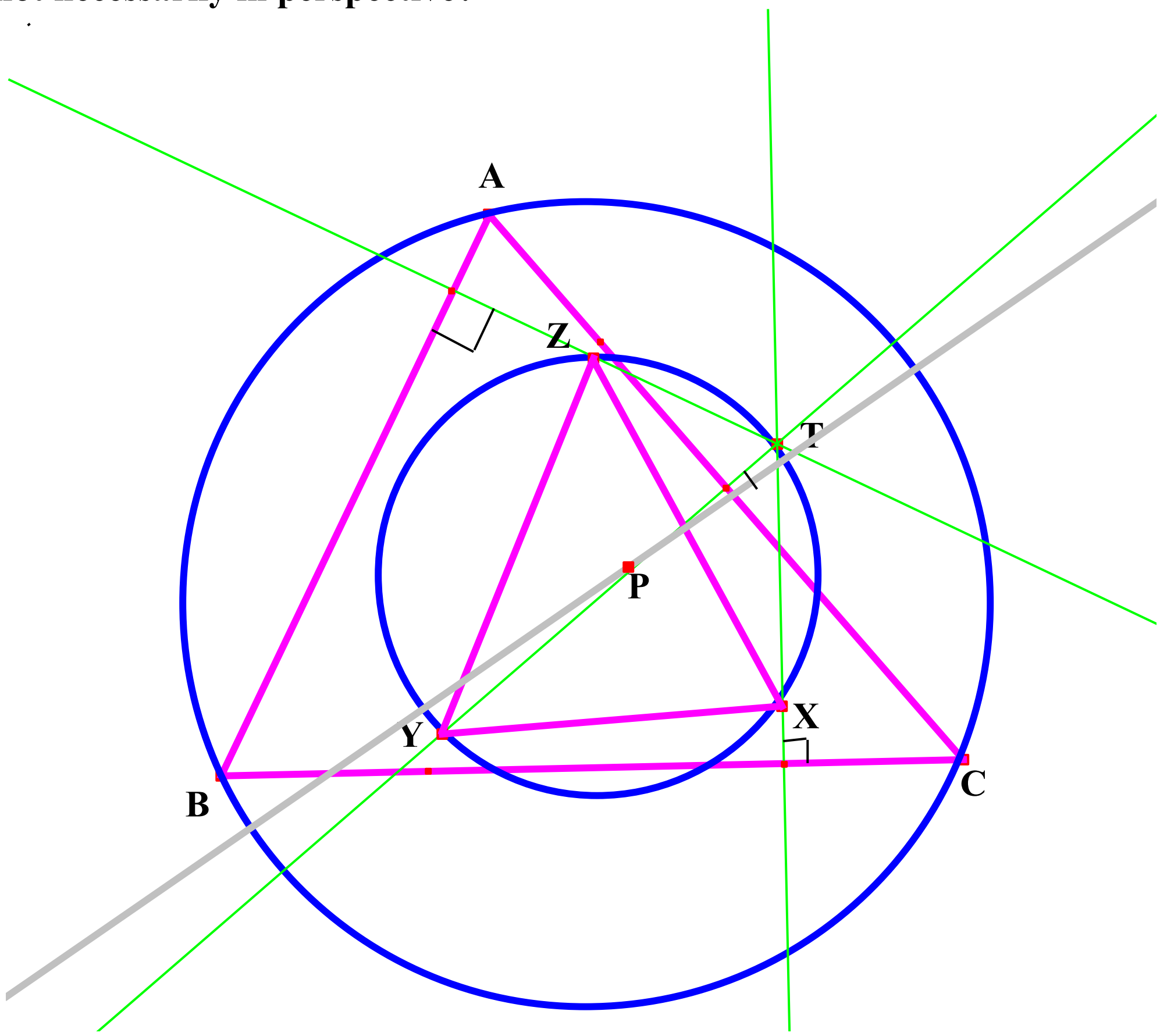

**Fig. 7**
**Constructing an indirect similarity**

We first provide a construction. Take a triangle *ABC*, draw any circle (other than the circumcircle) and choose any point *T* lying on it. Draw lines through *T* parallel to the altitudes of triangle *ABC* to meet the circle again at points *X*, *Y*, *Z*. See Fig. 7. Simple angle arguments show that triangle *XYZ* is indirectly similar to triangle *ABC*. As *T* moves around this circle the triangle *XYZ* also rotates around the circle and CABRI shows there may be up to three positions of *T* for which triangles *ABC* and *XYZ* are also in perspective, but that is not a feature that interests us in this paragraph. The construction is such that the lines through *X*, *Y*, *Z* perpendicular to *BC*, *CA*, *AB* respectively are concurrent at *T* and hence triangles *XYZ* and *ABC* are not only indirectly similar, they are orthologic.

**Theorem 7.1**
Every pair of indirectly similar triangles are orthologic and the orthologic centres lie on the circumcircles of the triangles.

*Proof*
Given a triangle *XYZ* that is indirectly similar to a triangle *ABC*, we can construct the circumcircle of *XYZ* and choose any point *T* on this circumcircle. We can then generate a triangle *X'Y'Z'* by drawing lines through *T* parallel to the altitudes of triangle *ABC*. Now all indirectly similar triangles inscribed in a given circle are congruent, so we can rotate *T* until triangles *X'Y'Z'* and *XYZ* coincide. The construction given at the outset of Section 7 now ensures that triangle *XYZ* is orthologic to *ABC* with orthologic centre at the final position of *T*.

Triangles *DEF* and *UVW* are now defined as in the Hagge construction and they are similar as one is the image of the other under the indirect similarity. Paralogic centres *S* and *s* may be defined, even when there is no perspective. However, when there are no perspectives, there are no rectangular hyperbolas. But it is true in the more general case that the double lines of inverse symmetry are parallel to the angle bisectors of corresponding sides of triangle *ABC* and *XYZ*.

When there is no perspective there is no Desargues axis of perspective *LMN*. However Cabri confirms that the axis *U'V'W'P* exists, as defined at the end of Section 2. The lines *KS* and *Ts* pass through the centres of their respective circles and so may justifiably be called orthologic diameters

## 8. What happens when the centre of inverse symmetry lies at the orthocentre

When Speckman wrote his paper [1] on indirectly similar triangles in perspective he kept his account general and did not take account of particular cases. In other words his account assumes that the centre of inverse symmetry *P*, the perspector *Q* and the orthocentres of triangles *ABC* and *XYZ* are distinct points. The result is that he missed two special cases that are highly interesting. When the perspector *Q* coincides with *H*, the orthocentre of *ABC*, one gets triangles inscribed in circles through *H* and this gap was filled three years later by the ground breaking work of Hagge [2]. What has never been considered in the period from 1905 to the present day is what happens when the centre of inverse symmetry *P* lies at *H*. In fact it produces cases that are just as interesting as the Hagge circles that arise when *Q* lies at *H*.

When *P* lies at *H* there are an infinite number of rectangular hyperbolas passing through *A*, *B*, *C* and *H* depending on the direction of the lines of inverse symmetry. And sure enough it turns out that if you choose ***any*** axis through *H* and then use it to form a triangle *XYZ* that is inversely similar to *ABC*, *XYZ always turns out to be in perspective with ABC*. The circumcircle of triangle *XYZ* obviously cannot be reduced by dilation through *H* to become a non-degenerate Hagge circle. In other words there are a whole set of pairs of indirectly similar and in perspective triangles and a whole set of circles, whose properties have not been studied before. Of course, since Speckman just starts with a given set of points *A*, *B*, *C*, *P* and *Q*, his results are applicable when *P* lies at *H*, but he never considered what additional properties are true when either *P* or *Q* lies at *H*. And these additional properties are, as we shall see, quite substantial. Hagge circles are the outcome of what happens when *Q* lies at *H*. We now describe what happens when *P* lies at *H*.

First, if you reflect $ABC$ in any line through $H$, and use any enlargement factor $k$ ($k \neq 0,1$) for an enlargement (or reduction) centre $H$, then the resulting triangle $XYZ$ is always similar to and in perspective with triangle $ABC$, with vertex of perspective at some point $Q$. Now let $AH$, $BH$, $CH$ meet the circumcircle of $ABC$ at $D$, $E$, $F$ respectively. Through $D$ draw a line parallel to $AX$ to meet $XP$ at $U$ and define $V$, $W$ similarly. Then it transpires that points $U$, $V$, $W$ lie on the circumcircle of triangle $XYZ$ and that triangle $UVW$ is the image of triangle $DEF$ under the indirect similarity. We know from work in Section 7, since $XYZ$ is indirectly similar to $ABC$, there must exist points $T'$ and $T$ that are orthologic centres of the two triangles and that these points lie on the circles $XYZ$ and $ABC$ respectively and have the properties that the perpendiculars from $T'$ on to $BC$, $CA$, $AB$ pass through $X$, $Y$, $Z$ respectively and the perpendiculars from $T$ on to $YZ$, $ZX$, $XY$ pass through $A$, $B$, $C$ respectively. It also follows that $T'$ lies on the rectangular hyperbola $XYZHQ$ and $T$ lies on the rectangular hyperbola $ABCHQ$. See Fig. 8 for illustration of all these properties. The earlier properties involving triangles $XYZ$, $DEF$ and $UVW$ require analytic proof and we now turn our attention to setting up this analysis.

**Choosing co-ordinates**

In dealing with a triangle it is possible to choose three co-ordinates, such as the angular dispositions of the vertices and thereby make use of symmetry, or it is possible, by choosing the axes in some preferred directions, to choose just two co-ordinates. In the latter case symmetry is lost, but the expressions for the equations of lines and co-ordinates of points may be less complicated. In the present calculation which features an axis through $H$ and dilation through $H$ it is desirable for this point to be taken as origin. We therefore present a calculation in which the vertices have co-ordinates as follows:

$A(-2-2vw, 0)$, $B(-2vw, 2v)$, $C(-2vw, 2w)$. It may now be checked that $H$ is the origin and the circumcentre $O$ has co-ordinates $(-1-3vw, v+w)$. The radius of the circumcircle $R$ is given by $R^2 = (1+v^2)(1+w^2)$.

The equations of $BC$, $CA$, $AB$ are respective

$$x + 2vw = 0, \qquad (8.1)$$
$$y - wx = 0, \qquad (8.2)$$
$$y - vx = 2v(1 + vw), \qquad (8.3)$$

and the equations of $AH$, $BH$, $CH$ are respectively

$$y = 0, \qquad (8.4)$$
$$wy + x = 0, \qquad (8.5)$$
$$vy + x = 0. \qquad (8.6)$$

The equation of the circumcircle $\Gamma$ is

$$x^2 + y^2 + 2(1 + 3vw)x - 2(v + w)y + 8vw(1 + vw) = 0. \qquad (8.7)$$

**The indirect similarity**

First perform the dilation through $H$ by a factor of $k$ ($k \neq 0, 1$) and then the images of $A$, $B$, $C$ are respectively $X'$, $Y'$, $Z'$ with co-ordinates $X'(-2k(1 + vw), 0)$, $Y'(-2kvw, 2kv)$, $Z'(-2kvw, 2kw)$. Since $H$ is the origin we may suppose the axis of reflection through $H$ has equation $y = mx$. We quote the result that the image of the point with co-ordinates $(c, d)$ after reflection in this line has co-ordinates $\{1/(1 + m^2)\}(2dm + c(1 - m^2), 2cm - d(1 - m^2))$. Using the co-ordinates of $X'$, $Y'$, $Z'$ as $c$, $d$ we obtain the co-ordinates of the images $X$, $Y$, $Z$ as

$$X: \{1/(1 + m^2)\}(-2k(1 - m^2)(1 + vw), -4km(1 + vw)),$$
$$Y: \{2kv/(1 + m^2)\}(2m - w(1 - m^2), -2mw - (1 - m^2)),$$
$$Z: \{2kw/(1 + m^2)\}(2m - v(1 - m^2), -2mv - (1 - m^2)).$$

**The perspective**

We now prove that triangles $ABC$ and $XYZ$ are in perspective. Since this is true for all values of $k \neq 0, 1$ and for all values of $m$ this is a very general result that implies that any axis through $H$ may be used for the reflection and any enlargement factor ($k = 1$would mean that $AX$, $BY$, $CZ$ are parallel). An interesting case that we do not provide separate analysis for is when the axis through $H$ is the Euler line and the enlargement factor is 0.5. The image of the circumcircle is then the nine-point circle and the points $X$, $Y$, $Z$ are additional natural points on the nine-point circle.

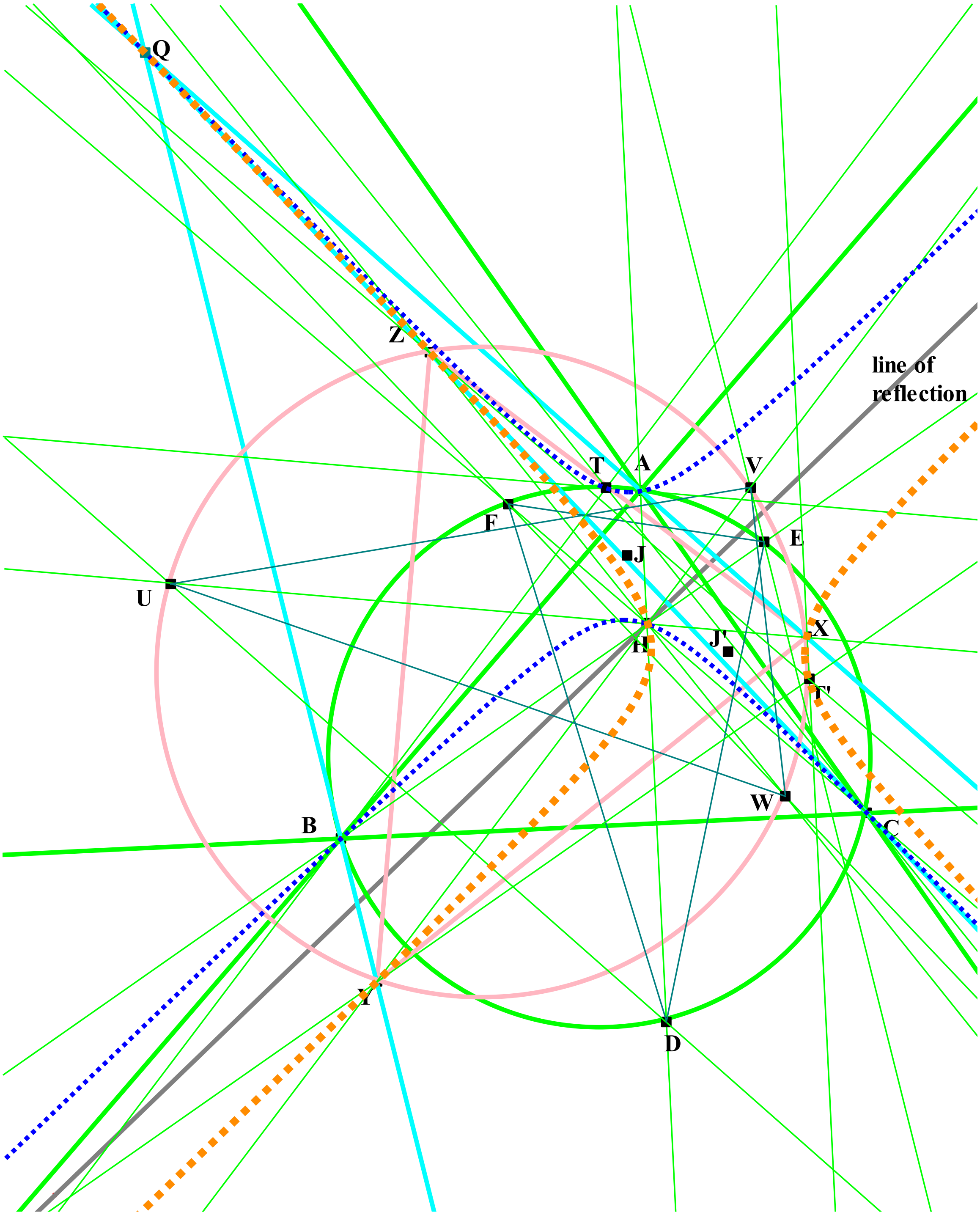

**Fig. 8**
**Indirectly similar triangles with reflection axis through the orthocentre**

The equations of $AX$, $BY$, $CZ$ turn out to be

$$2kmx + (1 + m^2 - k(1 - m^2))y + 4km(1 + vw) = 0, \qquad (8.8)$$

$$(1 + m^2 + k(1 + 2wm - m^2))x + (w(1 + m^2) - k(w - 2m - wm^2)y = 4kv(m(1 - w^2) - w(1 - m^2), \qquad (8.9)$$

$$(1 + m^2 + k(1 + 2vm - m^2))x + (v(1 + m^2) - k(v - 2m - vm^2)y = 4kw(m(1 - v^2) - v(1 - m^2), \qquad (8.10)$$

respectively.

These lines are concurrent at the point $Q$ with co-ordinates

$$Q\colon (4k/\{(1 - k^2)(1 + m^2)^2\})((k(m^4vw + m^3(v + w) + 2m^2 - m(v + w) + vw) + m^4vw + m^3(v + w) + m(v + w) - vw), - m(k(m^2(vw - 1) + 2m(v + w) - vw + 1) + (1 + m^2)(1 + vw))). \quad (8.11)$$

Of significance is the equation of circle $XYZ$ and this is

$$(1 + m^2)(x^2 + y^2) - 2k(m^2(3vw + 1) + 2m(v + w) - (1 + 3vw))x - 2k(m^2(v + w) - 2m(3vw + 1) - (v + w))y + 8k^2vw(1 + vw)(1 + m^2) = 0. \qquad (8.12)$$

**Triangle *UVW* and its properties**

We define points $U$, $V$, $W$ as the points of intersection of the lines $XH$, $YH$, $ZH$ with circle $XYZ$. This means that $UVW$ is inscribed in the same circle as triangle $XYZ$ and is in perspective with it. Since $X$, $Y$, $Z$ are the images of $A$, $B$, $C$ in the indirect similarity and $H$ is the fixed point of the similarity, it follows that if $AH$, $BH$, $CH$ meet the circumcircle at points $D$, $E$, $F$ respectively, then $U$, $V$, $W$ are the images of $D$, $E$, $F$ in the similarity and consequently triangles $DEF$ and $UVW$ are similar.

The equation of $XH$ is

$$(1 - m^2)y = 2mx. \qquad (8.13)$$

This meets circle $XYZ$ again at $U$, with co-ordinates $\{4k/(1 + m^2)\}(-(1 - m^2), -2m)$. The equation of $YH$ is

$$(m^2w + 2m - w)y = (m^2 - 2mw - 1)x. \qquad (8.14)$$

This meets circle $XYZ$ again at $V$, with co-ordinates

$$\{4kw(1 + vw)/(1 + m^2)(1 + w^2)\}(m^2w + 2m - w, m^2 - 2mw - 1).$$

The equation of $ZH$ is

$$(m^2v + 2m - v)y = (m^2 - 2mv - 1)x. \qquad (8.15)$$

This meets circle $XYZ$ again at $W$, with co-ordinates

$$\{4kv(1 + vw)/(1 + m^2)(1 + v^2)\}(m^2v + 2m - v, m^2 - 2mv - 1).$$

The equation of $AH$ is

$$y = 0. \tag{8.16}$$

This meets the circumcircle at the point $D$, with co-ordinates $(-4vw, 0)$.

The equation of $BH$ is

$$wy + x = 0. \tag{8.17}$$

This meets the circumcircle at the point $E$, with co-ordinates

$$\{4w(1 + vw)/(1 + w^2)\}(-w, 1).$$

The equation of $CH$ is

$$vy + x = 0. \tag{8.18}$$

This meets the circumcircle at the point $F$, with co-ordinates

$$\{4v(1 + vw)/(1 + w^2)\}(-v, 1).$$

Having defined $U$, $V$, $W$ in the above manner their main property is that the line through $D$ parallel to $AX$ passes through $U$, the line through $E$ parallel to $BY$ passes through $V$ and the line though $F$ parallel to $CZ$ passes through $W$.

The line $AX$ has equation given by (8.8), so the line parallel to $AX$ through $D$ has equation

$$2kmx + (1 + m^2 - k(1 - m^2))y + 8kmvw = 0. \tag{8.19}$$

It is now straightforward to show this meets the line $XH$, with equation (8.13), at the point $U$ with co-ordinates $\{4k/(1 + m^2)\}(-(1 - m^2), -2m)$. Similarly the line parallel to $BY$ through $E$ meets $YH$ at $V$ and the line parallel to $CZ$ through $F$ meets $ZH$ at $W$.

Flat 4,
Terrill Court,
12-14 Apsley Road,
BRISTOL BS8 2SP